\input amstex
\magnification=\magstep1 
\baselineskip=13pt
\documentstyle{amsppt}
\vsize=8.7truein \CenteredTagsOnSplits \NoRunningHeads
\def\conv{\operatorname{conv}}

\def\BB{\Cal{B}}

\def\L{\Cal{L}}

\topmatter

\title Centrally symmetric polytopes with many faces
\endtitle
\author Alexander Barvinok, Seung Jin Lee, and Isabella Novik \endauthor
\address Department of Mathematics, University of Michigan, Ann Arbor,
MI 48109-1043  \endaddress
\email barvinok$\@$umich.edu \endemail
\address Department of Mathematics, University of Michigan, Ann Arbor,
MI 48109-1043 \endaddress
\email lsjin$\@$umich.edu \endemail
\address Department of Mathematics, University of Washington, Seattle, WA 98195-4350 \endaddress
\email novik$\@$math.washington.edu \endemail 
\date November 2011 \enddate
\thanks  The research of the first and second authors was partially supported by 
NSF Grant DMS-0856640; 
the research of the third author was partially supported by NSF Grants DMS-0801152 and DMS-1069298.
\endthanks 
\abstract We present explicit constructions of centrally symmetric polytopes with many faces:
(1) we construct a $d$-dimensional centrally symmetric polytope $P$ with about 
$3^{d/4} \approx (1.316)^d$ vertices such that 
every pair of non-antipodal vertices of $P$ spans an edge of $P$, (2) 
for an integer $k \geq 2$, we construct a $d$-dimensional centrally symmetric polytope $P$ of an arbitrarily high
dimension $d$ and
with an arbitrarily large number $N$ of vertices such that for some $0 <\delta_k<1$
at least $(1-(\delta_k)^d){N \choose k}$ $k$-subsets of the set of vertices span faces of $P$,
and (3) for an integer $k \geq 2$ and $\alpha>0$, 
we construct a centrally symmetric polytope $Q$ with an arbitrarily large 
number of vertices $N$ and of dimension $d=k^{1+o(1)}$ such that at least $\left(1-k^{-\alpha}\right) {N \choose k}$
$k$-subsets of the set of vertices span faces of $Q$.
\endabstract
\endtopmatter

\document

\head 1. Introduction and main results \endhead

A {\it polytope} is the convex hull of a set of finitely many points in ${\Bbb R}^d$. 
A polytope $P \subset {\Bbb R}^d$ is {\it centrally symmetric} if $P =-P$. 
We present explicit constructions of centrally symmetric polytopes with many faces.
Recall that a {\it face} of a convex body is the intersection of the body with a 
supporting affine hyperplane, see, for example, Chapter II of \cite{Ba02}.

A construction of {\it cyclic polytopes}, which goes back to Carath\'eodory \cite{Ca11} 
and was studied by Motzkin \cite{Mo57} and Gale \cite{Ga63}, presents a family of 
polytopes in ${\Bbb R}^d$ with an arbitrarily large number $N$ 
of vertices, such that the convex hull of every set of $k \leq d/2$ vertices 
is a face of $P$. Such a polytope is obtained 
as the convex hull of a set of $N$ distinct points on the moment curve 
$\left(t, t^2, \ldots, t^d\right)$ in ${\Bbb R}^d$.

The situation with centrally symmetric polytopes is far less understood. 
A centrally symmetric polytope $P$ is called $k$-{\it neighborly} 
if the convex hull of every set $\{v_1, \ldots, v_k\}$ of $k$ vertices of $P$, not containing a 
pair of antipodal vertices $v_i=-v_j$, is a face of $P$. 
In contrast with polytopes without symmetry, even 2-neighborly centrally symmetric 
polytopes cannot have too many vertices: it was shown in \cite{LN06} that no
$d$-dimensional 2-neighborly centrally symmetric polytope has more than $2^d$ vertices.
Moreover, as was verified in \cite{BN08}, 
the number $f_1(P)$ of edges (1-dimensional faces)
of an arbitrary centrally symmetric polytope $P \subset {\Bbb R}^d$ 
with $N$ vertices satisfies 
$$f_1(P) \ \leq \ {N^2 \over 2} \left(1-2^{-d}\right).$$
Let $f_k(P)$ denote the number of $k$-dimensional faces 
of a polytope $P$. Even more generally, \cite{BN08} proved that for a $d$-dimensional 
centrally symmetric polytope $P$ with $N$ vertices, 
$$f_{k-1}(P) \ \leq \ {N \over N-1} 
\left(1 - 2^{-d}\right) {N \choose k}, \quad \text{provided} \quad k \leq d/2.$$
In particular, as the number $N$ of vertices grows while the dimension $d$ 
of the polytope stays fixed, the fraction of $k$-tuples $v_1, \ldots, v_k$ 
of vertices of $P$ that do not form the vertex set of a $(k-1)$-dimensional face of $P$
remains bounded from below by roughly $2^{-d}$.

Besides being of intrinsic interest, centrally symmetric polytopes with many faces 
appear in problems of sparse signal reconstruction, see \cite{Do04}, \cite{RV05}, and also Section 5.
Typically, such polytopes are obtained through a randomized 
construction, for example, as the orthogonal projection of a high-dimensional 
cross-polytope (octahedron) onto a random subspace, see \cite{LN06} and \cite{DT09}. 

In this paper, we present explicit deterministic constructions. 
First, we construct a $d$-dimensional 2-neighborly centrally symmetric polytope 
with roughly $3^{d/4} \approx (1.316)^d$ vertices. Then, for any fixed $k \geq 2$, 
we verify (again by presenting an explicit construction)
 that there exists $0<\delta_k<1$ such that for an arbitrarily large $d$ and for an
arbitrarily large even $N$, there is
a $d$-dimensional centrally symmetric polytope $P$ with 
 $N$ vertices satisfying
$$f_{k-1}(P) \ \geq \ \left(1 - (\delta_k)^d \right) {N \choose k}. $$
Our construction guarantees that one can take $$\text{any} \quad \delta_2 > 3^{-1/4} \approx 0.77 \quad \text{and any} 
\quad \delta_k > \left(1-5^{-k+1}\right)^{5/(24k+4)}  
\text{ for  $k>2$}
$$
provided $N$ and $d$ are sufficiently large.
Finally, for an integer $k \geq 2$ and $\alpha>0$ we construct a centrally symmetric polytope $Q$ of dimension $k^{1+o(1)}$ 
with an arbitrarily large number of vertices $N$ such that 
$$f_{k-1}(Q)\ \geq \ \left(1-k^{-\alpha}\right) {N \choose k}.$$

We note that the random projection construction cannot produce polytopes with the last two properties
 since if $N$ is very large compared to $d$, the projection of a cross-polytope 
in ${\Bbb R}^N$ onto a random $d$-dimensional subspace
is very close to a Euclidean ball, and hence has few faces relative 
to the number of vertices, cf.~\cite{DT09}.
Our constructions are based on the symmetric moment curve introduced in \cite{BN08} and 
further studied in \cite{B+11}.

\subhead (1.1) The symmetric moment curve \endsubhead 
We define the {\it symmetric moment curve} $U_k(t) \in {\Bbb R}^{2k}$ by 
$$U_k(t)=\Bigl(\cos t,\ \sin t,\ \cos 3t,\ \sin3t,\ \ldots,\ 
\cos (2k-1)t,\ \sin(2k-1)t \Bigr) \tag1.1.1$$
for $t \in {\Bbb R}$. Since 
$$U_k(t) = U_k(t + 2\pi) \quad \text{for all} \quad t,$$
from this point on, we consider $U_k(t)$ to be defined on the unit circle 
$${\Bbb S} ={\Bbb R} /2 \pi {\Bbb Z}.$$
We note that $t$ and $t+\pi$ form a pair of antipodal points for all $t \in {\Bbb S}$ and that 
$$U_k(t+\pi)=-U_k(t) \quad \text{for all} \quad t \in {\Bbb S}.$$

First, we construct a 2-neighborly centrally symmetric polytope using the curve
$$U_3(t)=\Bigl(\cos t,\ \sin t,\ \cos 3t,\ \sin 3t,\ \cos 5t,\ \sin 5t\Bigr).$$

\proclaim{(1.2) Theorem} For a non-negative integer $m$, consider the map
$$\Psi_m: {\Bbb S} \longrightarrow {\Bbb R}^{6(m+1)}  \quad \text{defined by} \quad 
\Psi_m(t)=\Bigl(U_3(t),\ U_3(3t), \ldots, U_3\left(3^m t\right)\Bigr).$$
Let $A_m \subset {\Bbb S}$ be the set of $4 \cdot 3^{m+1}$
equally spaced points,
$$A_m=\left\{ {2\pi j \over 4 \cdot 3^{m+1}}, \quad j=0, \ldots, 4 \cdot 3^{m+1}-1\right\},$$ 
and let
$$P_m=\conv\Bigl( \Psi_m(t): \quad t \in A_m \Bigr).$$
Then $P_m$ is a centrally symmetric polytope of dimension $d=4m+6$ that has
$4 \cdot 3^{m+1}$ vertices: $\Psi_m(t)$ for  $t \in A_m$.
Moreover, for $t_1, t_2 \in A_m$ such that $t_1\ne t_2$ and
$t_1 \ne t_2 +\pi \mod 2\pi$, the interval 
$$\bigl[ \Psi_m\left(t_1\right),\ \Psi_m\left(t_2\right)\bigr]$$ is an edge of $P_m$.
\endproclaim 

Our construction of a centrally symmetric polytope with $N$ vertices
and about $(1-3^{-d/4}){N \choose 2}$ edges for an arbitrarily large $N$ is 
a slight modification of the construction presented in Theorem 1.2 --- see Remark 3.2.
On the other hand, to construct a centrally symmetric polytope with many 
$(k-1)$-dimensional faces for $k>2$, we need to use the curve (1.1.1) to the full extent.

\proclaim{(1.3) Theorem} 
Fix an integer $k \geq 1$. For a non-negative integer $m$, consider the map
$\Psi_{k,m}: {\Bbb S} \longrightarrow {\Bbb R}^{6k(m+1)}$ defined by 
$$\Psi_{k,m}(t)=\Bigl(U_{3k}(t),\ U_{3k}(5t),\ \ldots,\ U_{3k}\left(5^m t\right)\Bigr).$$ 
For a positive even integer $n$, 
let $A_{m,n} \subset {\Bbb S}$ be the set of $n 5^m$ equally spaced points, 
$$A_{m,n}=\left\{ {2 \pi j \over n 5^m}:\quad j=0, \ldots, n5^m-1 \right\},$$
and let 
$$P=P_{k,m,n}=\conv\Bigl( \Psi_{k,m}(t): \quad t \in A_{m,n} \Bigr).$$
Then 
\roster
\item 
The polytope $P \subset {\Bbb R}^{6k(m+1)}$ is a centrally symmetric polytope 
 with $n 5^m$ distinct vertices:
$$ \Psi_{k,m}(t) \quad \text{for } t \in A_{m,n}$$
and of dimension $d\leq 6k(m+1)-2m\lfloor (3k+2)/5\rfloor$; moreover, if $n> 2(6k-1)$, then the dimension of $P$ is equal
to $6k(m+1)-2m\lfloor (3k+2)/5\rfloor$.
\item Let $t_1, \ldots, t_k $ be points chosen independently at random 
from the uniform distribution in $A_{m,n}$ (in particular, some of $t_i$ may coincide). 
Then the probability that
$$\conv\Bigl(\Psi_{k,m}\left(t_1\right), \ldots, \Psi_{k,m}\left(t_k\right) \Bigr)$$
is not a face of $P$ does not exceed 
$$\left(1-5^{-k+1}\right)^m.$$
\endroster
\endproclaim
We obtain the following corollary.

\proclaim{(1.4) Corollary} 
Let $P_{k,m,n}$ be the polytope of Theorem 1.3 with $N=n5^m$ vertices
and dimension $d \leq 6k(m+1)-2m\lfloor (3k+2)/5\rfloor$.
Then
$$f_{k-1}\left(P_{k,m,n}\right) \ \geq \ 
{N \choose k} -\left(1-5^{-k+1}\right)^m {N^k \over k!}.$$ 
\endproclaim

The construction of Theorem 1.3 produces a family of centrally symmetric 
polytopes of an increasing dimension $d$ and with an arbitrarily 
large number of vertices such that for any 
fixed $k \geq 1$, the probability $p_{d,k}$ that $k$ randomly chosen vertices 
of the polytope do not span a face decreases exponentially in $d$. 
However, it does not start doing so very quickly:
for instance, to make $p_{d,k} < 1/2$ we need to choose $d$ as high as $2^{\Omega(k)}$.

 Using a trick which the authors learned from Imre B\'ar\'any 
(cf.~Section 7.3 of  \cite{BN08}), we construct new families of 
polytopes with many faces of a reasonably high dimension. 
Namely, we can make $p_{d,k} <d^{-\alpha}$ for any fixed
 $\alpha>0$ by using $d$ as low as $k^{1+o(1)}$.

\proclaim{(1.5) Theorem} Fix positive integers $k, m, n$ and $r$, where $n$ is even.
Let $P=P_{k,m,n}$ be the polytope of Theorem 1.3, 
so that $P \subset {\Bbb R}^{6k(m+1)}$ is a centrally symmetric polytope 
with $n 5^m$ vertices. For $d=6kr(m+1)$, identify 
${\Bbb R}^d$ with a direct sum of $r$ copies of ${\Bbb R}^{6k(m+1)}$, 
each containing a copy of $P$. Let $Q$ be the convex hull of the $r$ copies of $P$;
in particular, $Q \subset {\Bbb R}^d$ is a centrally symmetric polytope with 
$r n 5^m$ vertices.

If
$$r \ < \ \min \left\{ (k+1)! , \ \left({5^{k-1} \over 5^{k-1}-1} \right)^m \right\},$$
then the probability that $r$ vertices of $Q$, chosen independently at random from 
the uniform distribution on the set set of vertices of $Q$, span a face of $Q$ is at least
$$\left(1 -{r \over (k+1)!}\right) \left(1- r \left(1-5^{-k+1}\right)^m \right).$$
\endproclaim 

If we now fix an $\alpha >0$ and choose in Theorem 1.5
$$k =\left\lceil {\beta \ln r \over \ln \ln r} \right\rceil 
\quad \text{and} \quad m=\left\lceil \beta 5^k \ln r \right\rceil,$$
then for a suitable $\beta =\beta(\alpha)>0$
 we obtain a centrally symmetric polytope $Q$ of dimension $r^{1+o(1)}$ and 
with an arbitrarily large number $N$ of vertices such that $r$ random vertices of $Q$
 span a face of $Q$ with probability at least $1-r^{-\alpha}$. As in Corollary 1.4, 
we have $f_{r-1}(Q) \geq \left(1 -r^{-\alpha}\right){N \choose r}$.

In Section 2, we summarize the properties of the symmetric moment curve (1.1.1) 
and review several basic combinatorial facts needed for our proofs.
We then prove Theorem 1.2 in Section 3 and Theorems 1.3 and 1.5 in Section 4. 
In Section 5, we sketch connections to error-correcting codes.

\head 2. Preliminaries \endhead 

We utilize the following result of \cite{B+11} concerning 
the symmetric moment curve~(1.1.1). 
\proclaim{(2.1) Theorem} Let $\BB_k \subset {\Bbb R}^{2k}$,
$$\BB_k=\conv\Bigl(U_k(t): \quad t \in {\Bbb S}\Bigr),$$
be the convex hull of the symmetric moment curve. Then
for every positive integer $k$ there exists a number 
$${\pi \over 2} \ < \ \alpha_k \ < \ \pi$$
such that for an arbitrary open arc $\Gamma \subset {\Bbb S}$ of length $\alpha_k$ 
and arbitrary distinct $n \leq k$ points $t_1, \ldots, t_n \in \Gamma$, the set 
$$\conv\Bigl(U_k\left(t_1\right), \ldots, U_k\left(t_n\right)\Bigr)$$
is a face of $\BB_k$.
\endproclaim

\noindent For $k=2$ with $\alpha_2=2\pi/3$ this result is due to Smilansky \cite{Sm85}.

We will also need the following technical lemma.
\proclaim{(2.2) Lemma}
Let $t_1, \ldots, t_{2k} \in {\Bbb S}$ be distinct points no two of which are antipodal. 
Then the set of vectors
$$\{U_k\left(t_1\right), \ldots, U_k\left(t_{2k}\right)\}$$
is linearly independent.
\endproclaim 
\demo{Proof} 
Seeking a contradiction, we assume that these $2k$ vectors are linearly dependent. Then
they span a proper subspace in ${\Bbb R}^{2k}$, and hence 
there is a non-zero vector $C \in {\Bbb R}^{2k}$ that is orthogonal to 
all these vectors.

Consider the following trigonometric polynomial 
$$f(t)=\big\langle C,\ U_k(t) \bigr\rangle,$$
where $\langle \cdot, \cdot \rangle$ is the standard scalar product in ${\Bbb R}^{2k}$.
Then $f(t) \not\equiv 0$ and $t_1, \ldots, t_{2k}$ are distinct roots of $f(t)$. 
Since $f(t+\pi)=-f(t)$, we conclude that $f(t)$ has at least $4k$ roots on 
the circle ${\Bbb S}$. On the other hand, substituting $z=e^{it}$, we can write 
$$f(t)={p(z)\over z^{2k-1}},$$
where $p$ is a polynomial with $\deg p \leq 4k-2$, see \cite{BN08} and \cite{B+11}. 
Hence $p(z)$ has at least $4k$ distinct 
roots on the circle $|z|=1$ and we must have $p(z) \equiv 0$, which is a contradiction.
{\hfill \hfill \hfill}\qed
\enddemo

We will also be using the following two well-known facts.

 First, if $P$ is a polytope and $F$ is a face of $P$, then $F$
is a polytope: it is the convex hull of the vertices of $P$ that lie in $F$.
Moreover, every face of $F$ is also a face of $P$.

Second, if $T: {\Bbb R}^d \longrightarrow {\Bbb R}^k$ is a linear transformation and 
$P \subset {\Bbb R}^d$ is a polytope, 
then $Q=T(P)$ is a polytope and for every face $F$ of $Q$ the inverse image of $F$, 
$$T^{-1}(F)=\Bigl\{x \in P: \quad T(x) \in F \Bigr\},$$
is a face of $P$.
This face is the convex hull of the vertices of $P$ mapped by $T$ into vertices of $F$.

Finally, to estimate the dimension of the polytope $P_{k,n,m}$ in Theorem 1.3 we will 
rely on the following combinatorial lemma. For a set $U$ of integers and a constant
$c$, we define $cU:=\{cu \ : \ u\in U\}$.

\proclaim{(2.3) Lemma}
Let $K$ be the set of all odd integers in the closed interval $[1, 6k-1]$, and let 
$$T=\bigcup_{j=0}^{m} 5^{j}K.$$
Then 
$$ |T| = 3k(m+1)-m\lfloor (3k+2)/5\rfloor.$$
\endproclaim 
\demo{Proof} 
Denote by $X$ the set of all elements of $K$ that are not divisible by 5, and by $S$ the complement of $X$ in $K$.
Then the sets  $X, \ 5X,\  5^2X,\ \cdots,\  5^m X$ are pairwise disjoint and their union
consists of all elements of $T$ that are not divisible by $5^{m+1}$. On the other hand, every element of 
$T$ that is divisible by $5^{m+1}$ is of the form $5^m s$ for some $s\in S$ and every element of the 
form $5^ms$ for $s \in S$ belongs to $T$ and is divisible by $5^{m+1}$. Thus
$$T=\left(\bigcup_{j=0}^m 5^jX\right) \cup 5^mS,$$
 and the sets in the above union are pairwise disjoint. Hence
$$|T|=(m+1)|X|+|S|=(m+1)|K|-m|S|.$$
The statement now follows from the fact that there are $3k$ elements in $K$ and that exactly 
$\lfloor (3k+2)/5\rfloor$ of them are divisible by 5. {\hfill \hfill \hfill}\qed
\enddemo

\head 3. Centrally symmetric 2-neighborly polytopes \endhead 

\subhead (3.1) Proof of Theorem 1.2 \endsubhead
The transformation 
$$t \longmapsto t+ \pi \mod 2\pi$$
maps the set $A_m$ onto itself. 
Since $\Psi_m(t+\pi)=-\Psi_m(t)$, the polytope $P_m$ is centrally symmetric.
 Consider the projection ${\Bbb R}^{6(m+1)} \longrightarrow {\Bbb R}^6$ 
that forgets all but the first $6$ coordinates.
Then the image of $P_m$ is the polytope 
$$Q_m=\conv\Bigl( U_3(t): \quad t \in A_m \Bigr). \tag3.1.1$$
By Theorem 2.1, the 
polytope $Q_m$ has $4 \cdot 3^{m+1}$ distinct vertices: $U_3(t)$ for $t \in A_m$. Furthermore,
the inverse image of each vertex $U_3(t)$ of $Q_m$ in $P_m$
consists of a single vertex $\Psi_m(t)$ of $P_m$.
Therefore, $\Psi_m(t)$ for $t \in A_m$ are all the vertices of $P_m$ without duplicates.

To compute the dimension $d$ of $P_m$, we observe that for all $t \in {\Bbb S}$,
the third coordinate of $U_3(t)$ coincides with the first coordinate of $U_3(3t)$ 
while the fourth coordinate of $U_3(t)$ coincides with the second coordinate of $U_3(3t)$. 
Therefore, the polytope $P_m$ lies in a subspace, denote it by $\L$, of codimension $2m$, and hence $\dim P_m\leq 4m+6$. 
If the dimension of $P_m$ is strictly smaller than $4m+6$, then $P_m$ lies in an affine hyperplane of $\L$. 
As in the proof of Lemma 2.2, such an affine hyperplane corresponds 
to a trigonometric polynomial $f(t)$ of degree $5\cdot 3^m$
that has at least $4\cdot 3^{m+1}= 12\cdot 3^m$ roots (all points of $A_m$). This is however impossible, as no 
nonzero trigonometric polynomial of degree  $5\cdot 3^m$ has more than
$$
2\cdot 5\cdot 3^m = 10\cdot  3^m <  12\cdot 3^m$$
roots (cf.~the proof of Lemma 2.2). We conclude that $\dim P_m=4m+6$.

We prove that $P_m$ is 2-neighborly by induction on $m$. 
It follows from Lemma~2.2 that $P_0$ is the convex hull of 
a set consisting of six linearly independent vectors and their opposite vectors. 
Combinatorially, $P_0$ is a 6-dimensional
cross-polytope and hence the induction base is established. 

Suppose now that $m \geq 1$. Let $t_1, t_2 \in A_m$ be such that  
$$t_1 \ \ne \ t_2,\ t_2+\pi \mod 2 \pi.$$
Then there are two cases to consider.
\bigskip
{\sl Case I:} 
$\  t_1 - t_2 \in(-{\pi \over 2}, {\pi \over 2}) \mod 2 \pi$,
\medskip
\noindent{and}
\medskip
{\sl Case II:}  
$\  t_1 -t_2 \in(-\pi,-{\pi \over 2}] \cup [{\pi \over 2}, \pi) \mod 2 \pi .$
\bigskip
In the first case, consider the polytope $Q_m$ defined by (3.1.1) and the projection 
$P_m \longrightarrow Q_m$ as above. By Theorem 2.1,
$$\bigl[ U_3\left(t_1\right), \ U_3\left(t_2\right)\bigr]$$
is an edge of $Q_m$. Since the inverse image of a vertex $U_3(t)$ of $Q_m$ in $P_m$
 consists of a single vertex $\Psi_m(t)$ of 
$P_m$, we conclude that  
$$\bigl[ \Psi_m\left(t_1\right), \ \Psi_m\left(t_2\right) \bigr]$$
is an edge of $P_m$.

In the second case, consider the map $\phi: A_m \longrightarrow A_{m-1}$,
$$\phi(t) =3t \mod 2\pi.$$
Then 
$$\phi\left(A_m\right)=A_{m-1}$$
and for every $t$ the inverse image of $t$, $\phi^{-1}(t)$, consists of $3$ 
equally spaced points from $A_m$. 
In addition, we have 
$$\phi\left(t_1\right) \ne \phi\left(t_2\right) + \pi \mod 2\pi,$$
although we may have $\phi\left(t_1\right)=\phi\left(t_2\right)$. 
In any case, by the induction hypothesis, the interval
(possibly contracting to a point)
$$\bigl[ \Psi_{m-1}\left(3 t_1\right),\ \Psi_{m-1}\left(3 t_2\right) \big] \tag3.1.2$$
is a face of $P_{m-1}$. 

Let us consider the projection ${\Bbb R}^{6(m+1)} \longrightarrow {\Bbb R}^{6m}$ 
that forgets the first 6 coordinates.
The image of $P_m$ under this projection is $P_{m-1}$, 
and since (3.1.2) is a face of $P_{m-1}$, the set 
$$\aligned \conv\Bigl( \Psi_m \left(x_{ij}\right): \quad 
\phi\left(x_{ij}\right)=\phi\left(t_i\right) \quad 
\text{for} \quad &i=1, 2 \\ \text{and} \quad &j=1,2,3\Bigr) \endaligned\tag3.1.3$$
is a face of $P_m$ (it is the inverse image of (3.1.2) under this projection). 
However, the face (3.1.3) is a convex 
hull of at most six distinct points no two of which are antipodal. Since by Lemma 2.2,
any set of at most six distinct points $U_3\left(x_{ij}\right)$ no two of which are
antipodal is linearly independent, 
the face (3.1.3) is a simplex. Therefore,  
$$\bigl[ \Psi_m\left(t_1\right), \ \Psi_m\left(t_2\right) \big]$$
is a face of (3.1.3), and hence of $P_m$.
{\hfill \hfill \hfill} \qed

\remark{(3.2) Remark} Tweaking the construction of Theorem 1.2, allows us to
produce $d$-dimensional centrally symmetric polytopes
with an arbitrarily large number $N$ of vertices that have at least
$\left(1-(\delta_2)^d\right) {N \choose 2}$ edges, where one can choose any
$\delta_2 > 3^{-1/4} \approx 0.77$ 
for all sufficiently large $N$ and $d$.

To do so, fix an integer $s\geq 3$, and consider the curve $\Psi_m$ 
as in Theorem 1.2. However, instead of
working with the set $A_m$ as in the proof Theorem 1.2,  start with the
set $$W_0=\left\{\frac{\pi j}{2}: \quad j=0,1,2,3\right\}$$
 of $4$ equally spaced points on ${\Bbb S}$. Now replace each point $t$ of
$W_0$ by a cluster of $s$ points on ${\Bbb S}$ that lie very close to $t$.
Moreover, do it in such a way, that the resulting subset of ${\Bbb S}$,
which we denote by $W_0^s$, is centrally symmetric. For $m\geq 1$, define
$W_m^s$ recursively by 
$$
W_{m}^{s} :=\phi^{-1}(W_{m-1}^{s}), \quad \text{where} \quad \phi(x)=3x \mod 2\pi.
$$
Thus $W_{m}^{s}$ consists of $4\cdot 3^m$ clusters of $s$ points each.

We claim that the polytope 
$$P_{m}^{s}:=\conv \left(\Psi_m(t):
\quad t\in W_{m}^{s}\right)$$ 
is a centrally symmetric polytope of dimension
$d=4m+6$, with $N=N(s)=4s \cdot 3^{m}$ vertices, and such that for every
two distinct points $t_1, t_2 \in W_{m}^{s}$, the interval 
$[\Psi_m(t_1), \Psi_m(t_2)]$ is an edge of $P_{m}^{s}$, 
provided $t_1$ and $t_2$ are not from antipodal clusters. 
The proof of this claim is identical to the proof of Theorem 1.2,
except that for the base case (the case of $m=0$) we appeal to Theorem 2.1.

Thus each vertex of $P_{m}^{s}$ is incident to all other vertices
except itself and (possibly) the $\Psi_m$-images of the $s$ points from the
antipodal cluster. Therefore, the polytope $P_{m}^{s}$ has at
least 
$$\frac{N(N-s-1)}{2}={N \choose 2}\left(1-\frac{s}{N-1}\right) \approx 
{N \choose 2}\left(1-\frac{1}{4\cdot 3^{m}}\right)$$ 
edges. Taking an arbitrarily large $s$ yields the promised result on 
$\delta_2$. {\hfill \hfill \hfill} \qed
\endremark

\head 4. Centrally symmetric polytopes with many faces \endhead

\subhead (4.1) Proof of Theorem 1.3 \endsubhead
We observe that the transformation 
$$t \longmapsto t +\pi \mod 2\pi$$
maps the set $A_{m, n}$ onto itself and that 
$$\Psi_{k, m}(t+\pi)=-\Psi_{k,m}(t) \quad \text{for all} \quad t \in {\Bbb S}.$$
Hence $P$ is centrally symmetric. 
Consider the projection ${\Bbb R}^{6k(m+1)} \longrightarrow {\Bbb R}^{6k}$ that 
forgets all but the first $6k$ coordinates. Then the image of $P_{k,m,n}$ is the polytope
$$Q_{k,m,n}=\conv\Bigl(U_{3k}(t): \quad t \in A_{m,n}\Bigr). \tag4.1.1$$
By Theorem 2.1, the polytope $Q_{k,m,n}$ 
has $n5^m$ distinct vertices: $U_{3k}(t)$ for $t \in A_{m,n}$. 
Furthermore, the inverse image of each vertex $U_{3k}(t)$ of $Q_{k,m,n}$ in $P_{k,m,n}$ 
consists of a single vertex $\Psi_{k,m}(t)$ of $P_{k,m,n}$. Therefore,
$\Psi_{k,m,n}(t)$ for $t \in A_{m,n}$ are all the vertices of $P_{k,m,n}$
without duplicates. 

To estimate the dimension of $P=P_{k, m,n}$, we observe that for all $t\in {\Bbb S}$, 
the fifth coordinate of $U_{3k}(t)$ coincides with the first coordinate of 
$U_{3k}(5t)$ while the sixth coordinate of $U_{3k}(t)$ coincides with the 
second coordinate of $U_{3k}(5t)$, etc. Taking into account all coincidences of coordinates,
we infer from Lemma 2.3 that the polytope $P$ lies in 
a subspace of dimension $6k(m+1)-2m\lfloor (3k+2)/5\rfloor$, and hence $\dim P\leq  6k(m+1)-2m\lfloor (3k+2)/5\rfloor$. 
Moreover, if $n>2(6k-1)$, then an argument identical to the one used in the proof of
Theorem 1.2 (by counting roots of trigonometric polynomials) shows that $\dim P =  6k(m+1)-2m\lfloor (3k+2)/5\rfloor$.

We prove Part (2) by induction on $m$. 
The statement trivially holds for $m=0$. Let us assume that $m \geq 1$ and 
consider the map $\phi: A_{m,n} \longrightarrow A_{m-1,n}$ defined by 
$$\phi(t)=5t \mod 2 \pi.$$
Then 
$$\phi\left(A_{m, n}\right)=A_{m-1, n}$$
and for every $t \in A_{m-1, n}$,
the inverse image of $t$, $\phi^{-1}(t)$, consists of 5 equally spaced points from $A_{m, n}$.
We note that if $t$ is a random point uniformly distributed in $A_{m,n}$, 
then $\phi(t)$ is uniformly distributed in $A_{m-1, n}$. 
The proof of the theorem will follow from 
the following two claims.
\bigskip
{\sl Claim I.} Let $t_1, \ldots, t_k \in A_{m,n}$ be arbitrary, not necessarily distinct, points. If 
$$\conv\Bigl(\Psi_{k,m-1}\left(5t_i\right), \quad i=1, \ldots, k\Bigr) \tag4.1.2$$
is a face of $P_{k,m-1,n}$ then 
$$\conv\Bigl(\Psi_{k,m}\left(t_i\right), \quad i=1, \ldots, k \Bigr) \tag4.1.3$$
is a face of $P_{k, m, n}$.
\medskip
{\sl Claim II.} Let $s_1, \ldots, s_k \in A_{m-1,n}$ be arbitrary, not necessarily distinct, points. 
Then the conditional probability that
$$\conv\Bigl(\Psi_{k,m}(t_i): \quad i=1, \ldots, k\Bigr)$$
is not a face of $P_{k,m,n}$ given that 
$$\phi\left(t_i\right)=s_i \quad \text{for} \quad i=1, \ldots, k$$
does not exceed $1-5^{-k+1}$.
\bigskip
To prove Claim I, we consider the projection ${\Bbb R}^{6k(m+1)} \longrightarrow {\Bbb R}^{6km}$ that 
forgets the first $6k$ coordinates. 
The image of $P_{k,m,n}$ under this projection is $P_{k,m-1,n}$ and if 
(4.1.2) is a face of $P_{k,m-1,n}$ then 
$$\aligned \conv\Bigl(\Psi_{k,m}\left(x_{ij}\right):  \quad \phi\left(x_{ij}\right)=\phi\left(t_i \right)
\quad \text{for} \quad &i=1, \ldots, k \\ \text{and} \quad &j=1,2,3,4,5\Bigr)  \endaligned \tag4.1.4$$
is a face of $P_{k,m,n}$ as it is the inverse image of (4.1.2) under this projection. 
The face (4.1.4) is the convex hull of 
at most $5k$ distinct points and no two points $x_{ij}$ in (4.1.4) are antipodal.  
Since by Lemma 2.2 a set of up to $6k$  distinct points
$U_{3k}\left(x_{ij}\right)$ no two of which are antipodal is linearly independent, 
the face (4.1.4) is a simplex. Therefore, the set 
(4.1.3) is a face of (4.1.4), and hence also a face of $P_{k,m,n}$. 
Claim~I now follows.

To prove Claim II, we fix a sequence $s_1, \ldots, s_k \in A_{m-1,n}$ of not necessarily distinct 
points. Then there are exactly $5^k$ sequences $t_1, \ldots, t_k \in A_{m,n}$ of not necessarily 
distinct points such that $\phi\left(t_i\right)=s_i$ for $i=1, \ldots, k$. 
Choose an arbitrary $t_1$ subject to the condition $\phi\left(t_1\right)=s_1$.
Let $\Gamma \subset {\Bbb S}$ be a closed arc of length $2\pi/5$ centered at $t_1$. 
Then for $i=2, \ldots, k$ 
there is at least one $t_i \in \Gamma$ such that $\phi(t_i)=s_i$. By Theorem 2.1, 
for such a choice of $t_2, \ldots, t_k$, the set
$$\conv\Bigl( U_{3k}\left(t_i\right): \quad i=1, \ldots, k\Bigr) \tag4.1.5$$ 
is a face of the polytope $Q_{k,m,n}$ defined by (4.1.1). Considering the projection 
$$P_{k,m,n} \longrightarrow Q_{k,m,n}$$
as above, we conclude that (4.1.3) is a face of $P_{k,m,n}$ as it is the inverse image of~(4.1.5).

Hence the conditional probability that (4.1.3) is not a face is at most 
$${5^{k-1} -1 \over 5^{k-1}} = 1-5^{-k+1}.$$
{\hfill \hfill \hfill} \qed

\subhead (4.2) Proof of Corollary 1.4 \endsubhead
 Let us choose points $t_1, \ldots, t_k$ independently 
at random from the uniform distribution in $A_{m,n}$.
Then the probability that the points are all distinct is
$${(N-1) \cdots (N-k+1) \over N^{k-1}}.$$
From Theorem 1.3, the conditional probability that 
$$\conv\Bigl( \Psi_{k,m}\left(t_1\right), \ldots, \Psi_{k,m}\left(t_k\right) \Bigl) \tag4.3.1$$ 
is not a face, given that $t_1, \ldots, t_k$ are distinct, does not exceed 
$$\left(1 -5^{-k+1}\right)^m {N^{k-1} \over (N-1) \cdots (N-k+1)}.$$
Arguing as in the proof of Theorem 1.3 (Section 4.1), we conclude that if 
$t_1, \ldots, t_k$ are distinct and (4.3.1) is a face, then
that face is a $(k-1)$-dimensional simplex.
{\hfill \hfill \hfill} \qed

\subhead (4.3) Proof of Theorem 1.5 \endsubhead 
By construction, $Q$ is a centrally symmetric polytope whose vertex set consists of 
the vertices of the $r$ copies of $P$.
Let us pick $r$ vertices of $Q$ independently at random from 
the uniform distribution and let $k_i$ be the number 
of vertices picked from the $i$-th copy of $P$, $i=1, \ldots, r$. Then the 
probability that $k_i > k$ does not exceed 
$${r \choose k+1} r^{-k-1} \ < \ {1 \over (k+1)!}.$$
Therefore, the probability that $k_1, \ldots, k_r \leq k$ is at least $1-r/(k+1)!$. 
Now, the picked $r$ vertices span a face of 
$Q$ if and only if for all $i$ with $k_i>0$
the chosen $k_i$ vertices from the $i$-th copy of $P$ span a face of $P$.
The result then follows by Theorem 1.3.
{\hfill \hfill \hfill} \qed

\head 5. Connections to error-correcting codes \endhead

Here we briefly touch upon a well-known connection 
between centrally symmetric polytopes with many faces and the coding theory, 
see, for example, \cite{RV05}.

 Let ${\Bbb R}^N$ be $N$-dimensional Euclidean space with the standard 
basis $e_1, \ldots, e_N$ and the $\ell^1$-norm 
$$\|x\|_1=\sum_{i=1}^N \left|x_i\right| \quad \text{for} \quad x=\left(x_1, \ldots, x_N\right).$$
 Let $L \subset {\Bbb R}^N$ be a subspace, 
let $v_i$ be the orthogonal projection of $e_i$ onto $L$, and let 
$$P =\conv\Bigl( \pm v_i, \quad i=1, \ldots, N \Bigr)$$
be the orthogonal projection of the standard cross-polytope (octahedron) in ${\Bbb R}^N$
onto $L$.

Let $L^{\bot} \subset {\Bbb R}^N$ be the orthogonal complement of $L$. 
 Suppose that we are given a point $a \in {\Bbb R}^N$, $a=\left(a_1, \ldots, a_N\right)$,
which is obtained by changing (corrupting) some (unknown) $k$ 
coordinates of an unknown point $c \in L^{\bot}$, 
$c=\left(c_1, \ldots, c_N\right)$, and that our goal is to find $c$. 
One, by now standard, way of attempting to do that is to try 
to find $c$ as the solution to the linear programming 
problem of minimizing the function 
$$x \longmapsto \|x - a \|_1 \quad \text{for} \quad x \in L^{\bot}. \tag5.1$$
Indeed, let 
$$I_+ =\Bigl\{i: \quad c_i > a_i \Bigr\} \quad \text{and} \quad 
I_-=\Bigl\{i: \quad c_i < a_i \Bigr\}.$$
Then $c$ is the unique minimum point of (5.1) if 
$$\conv\Bigl( v_i \quad \text{for} \quad i \in I_+ \quad \text{and} 
\quad -v_i \quad \text{for} \quad i \in I_- \Bigr)$$
is a face of $P$. By constructing polytopes $P$ with many 
$(k-1)$-dimensional faces we produce subspaces $L^{\bot}$
with the property that the points of $L^{\bot}$ can be efficiently 
reconstructed from many of the different ways of corrupting some $k$ of their 
coordinates.

\enddocument